\documentclass[a4paper]{jpconf}
\usepackage{iopams}
\expandafter\let\csname equation*\endcsname\relax
\expandafter\let\csname endequation*\endcsname\relax

\usepackage{amsmath,empheq}
\usepackage{graphicx}
\usepackage{hyperref}
\usepackage{soul}
\usepackage[normalem]{ulem}
\usepackage{bbm}
\usepackage{cite}

	\newcommand{\ncd}{\newcommand}
	\ncd{\mrm} {\mathrm}
	\ncd{\beq} {\begin{equation}}
	\ncd{\eeq} {\end{equation}}

	\def\d{{\rm d}}

\begin{document}

	\title{Putting gravity in control}
	
		\author{C S Lopez-Monsalvo}
		\ead{cslopezmo@conacyt.mx}
		\address{Conacyt-Universidad Aut\'onoma Metropolitana Azcapotzalco\\Avenida San Pablo Xalpa 180, Azcapotzalco, Reynosa Tamaulipas, 02200 Ciudad de M\'exico, M\'exico}
		\author{I Lopez-Garcia, F Beltran-Carbajal and R Escarela-Perez}
		\address{Universidad Aut\'onoma Metropolitana Azcapotzalco\\Avenida San Pablo Xalpa 180, Azcapotzalco, Reynosa Tamaulipas, 02200 Ciudad de M\'exico, M\'exico}	
	
%		\author{I Lopez-Garcia}
%		\email{correo}
%		\affiliation{Universidad Aut\'onoma Metropolitana Azcapotzalco\\Avenida San Pablo Xalpa 180, Azcapotzalco, Reynosa Tamaulipas, 02200 Ciudad de México, M\'exico}

	\begin{abstract}
	  The aim of the present manuscript is to present a novel proposal in Geometric Control Theory inspired in the principles of  General Relativity and energy-shaping control.

	\end{abstract}

%\maketitle

It has been a pleasure to prepare this contribution to celebrate the two parallel lives of effort and inspiration in   Theoretical Physics of Rodolfo Gambini and Luis Herrera. Albeit different in scope, their research shares the common thread of gravity and its geometric principles, together with their ultimate consequences ranging from the very nature of space and time to the astrophysical implications of thermodynamics and gravity. It is thus our wish to take this opportunity to present a different perspective on the use of the same guiding principles that led to General Relativity and explore, even if very briefly, its implications within the area of control theory.

General relativity brought a new paradigm in the understanding of physical reality, establishing a deep connection between geometry and physics through the falsification of the motion of a test body due to a \emph{universal force} \cite{reichenbach2012philosophy} by free motion in a curved manifold. In such case, curvature becomes morally tantamount to a universal force field strength. This insight gave rise to a successful geometrization programme for field theories. However, due to the non-universal nature of the other known interactions, the geometrization is slightly different from that of gravity. 

In geometric field theories, problems are of two kinds: given the sources determine the curvature, or, given the curvature determine the trajectories of test particles. In this two type of problems, it is the source which determines the background geometry for the motion. Thus, in general, changing the source distribution changes the curvature of the spacetime manifold, and thus the motion of the test particles. Therefore, in principle, one could be able to reproduce any observed motion in space by means of a suitable distribution of energy. By placing appropriate sources here and there, one can \emph{control} the motion of test particles. 
Of course, one cannot simply engineer and energy-momentum tensor so that test particles follow our desired trajectories. Most of such energy distributions are indeed unphysical from the spacetime point of view. However, this very principle might be applied to a different physical setting, one for which the background geometry is not the spacetime manifold. Such is the case of the geometrization of classical mechanics and the control theory that can be built from it.

 We begin by briefly revisiting the geometrization  programme of classical mechanics. Let us consider a mechanical systems characterized by $n$ degrees of freedom (DoF) and defined by a Lagrangian function $L:T\mathcal{Q} \longrightarrow \mathbb{R}$. Here, the configuration space is an $n$ dimensional manifold $\mathcal{Q}$ whose tangent bundle is denoted by $T\mathcal{Q}$. The generic problem in geometric mechanics can be stated as follows. Given a pair of points $p_1, p_2 \in \mathcal{Q}$, find the curve $\gamma \subset \mathcal{Q}$
 	\beq
	\gamma:[0,1] \longrightarrow \mathcal{Q}, \quad \gamma(0)=p_1,\ \gamma(1)=p_2,
 	\eeq
for which the functional 
	\beq
	\label{action}
	S\left[ \gamma \right] = \int_0^1 L\left[\tilde\gamma(\tau)\right] \d\tau,
	\eeq
is an extremum. Here, $\tilde \gamma$ denotes canonical lift of $\gamma$ to $T\mathcal{Q}$.
	\beq
	\tilde \gamma:[0,1] \longrightarrow T\mathcal{Q}, \quad \tilde \gamma:\tau \longrightarrow \left[\gamma(\tau),\dot\gamma(\tau) \right].
	\eeq
The fact that the \emph{natural} evolution of a mechanical system connecting two given points of the configuration space $\mathcal{Q}$ follows precisely such a path is known as \emph{Hamilton's Principle}. Thus, natural motions solve the Euler-Lagrange equations
	\beq
	\label{EL1}
	\mathcal{E}(L) = 0,
	\eeq
 where $\mathcal{E}$ is the Euler operator \cite{olver2000applications}.
 
 We will consider systems  generated by Lagrangian functions of the form 
 	\beq
	L:T\mathcal{Q} \longrightarrow \mathbb{R}, \quad L = T - V,
	\eeq 
where $T$ and $V$ represent the kinetic and potential energies, respectively. Moreover, we will restrict our analysis to the case where the kinetic energy is defined in terms of a symmetric, non-degenerate and (usually) positive definite second rank tensor field
	\beq
	\label{KEMetric}
	M:T\mathcal{Q} \times T\mathcal{Q} \longrightarrow \mathbb{R}.
	\eeq
Namely, those for which the kinetic energy at a given point is written as
	\beq
	T\vert_p = \frac{1}{2} M \left(U_p,U_p\right),
	\eeq
where
	\beq
	\label{velocity}
	U_p = \left.\frac{\d \gamma}{\d\tau}\right\vert_{\tau_0}\in T_p\mathcal{Q}, \quad \text{with} \quad \gamma(\tau_0) = p \in \mathcal{Q}
	\eeq
is the velocity of the trajectory at the point $p$. Such systems are called \emph{natural} (cf. \cite{pettini2007geometry})

The tensor field \eqref{KEMetric} satisfies at each point of $\mathcal{Q}$ the properties of an inner product, promoting the configuration space into a Riemannian manifold. Thus, one could try to relate the curves extremizing the action functional $S$ to geodesics of $M$ in $\mathcal{Q}$. However, they only coincide in the case of \emph{free motion}, that is, when there is no potential energy nor external forces. Nonetheless, this observation gives rise to the problem of finding a metric tensor $G$ for the configuration space $\mathcal{Q}$ such that  the extremal curves of the action $S$ coincide with the geodesics of $G$ for a natural system defined by Lagrangian function with potential energy function $V$. 

Recalling that geodesics are themselves extrema of the arc-length functional of $\mathcal{Q}$, it is a straightforward exercise showing that (c.f. Chapter 4 in \cite{pettini2007geometry}), for conservative systems, the metric we are looking for is conformal to $M$. That is, the solutions to the geodesic equation
	\beq
	\label{geoG}
	\nabla_U U =  0,
	\eeq
where $\nabla$ is a connection compatible with the metric
	\beq
	\label{JacobiMetric}
	G = 2\left[E - V(p)\right] M,
	\eeq  
extremize the action functional \eqref{action}. Here, $E$ is the energy of the initial conditions and, by assumption, it is a constant of the motion. The metric \eqref{JacobiMetric} is called the \emph{Jacobi metric}. There is a difference, however, in the geometric origin of these curves. On the one hand, they are the paths followed by the system in a potential whilst, on the other, those are the free paths of the purely kinetic Lagrangian	
	\beq
	L_G = \frac{1}{2} G (U,U). %\quad \text{\bf check for square roots}
	\eeq

Now that we can identify the trajectories in configuration space of the \emph{natural motion} of a given system with geodesics in a Riemannian manifold, we would like to bend those paths so that the system evolves in a \emph{desired} manner, that is, we want to re-shape the geometry so that the \emph{desired evolution} corresponds to geodesic motion in a control Riemannian manifold. Thus, albeit both, classical mechanics and control theory  are based on the same dynamical principles, they do differ in their objectives and goals. On the one hand, the generic problem of classical mechanics is that of finding the integral curves to a given Lagrangian vector field whilst, on the other hand, in control theory one is interested in finding the control \emph{inputs}, generated by  properly located actuators, so that the integral curves of a given system follow a designed or desired path in configuration space. Moreover, in both cases, the problem of stability is of paramount relevance. In the following lines we present a proposal to address the stability problem in control theory from a Riemannian point of view. 
 
 The way one controls the evolution of a system is by directly acting upon a set of \emph{accesible} degrees of freedom (ADoF) \cite{lewis2004notes}
 	\beq
	\mathcal{D} = \left\{\hat e_{(i)}(p) \right\}_{i=1}^m, \quad \text{where} \quad \hat e_{(i)}(p) \in T_p \mathcal{Q} \quad \forall p \in \mathcal{Q},
	\eeq 
is the $i$th element of a frame over the configuration space so that we can place the control action directly into the Euler-Lagrange equations as
	\beq
	\label{controlEL}
	\mathcal{E}(L) = M^\sharp \left[u\right] \in T\mathcal{Q}, \quad u=\sum_{i=1}^m u_i \hat f^{(i)} \quad \text{with} \quad \hat f^{(i)} = M^\flat \left[\hat e_{(i)} \right],
	\eeq
where $u$ represents our control input as an applied external force acting on $\mathcal{D}$. Here, $M^\sharp$ and $M^\flat$ denote the musical isomorphisms between the tangent and co-tangent bundles defined by the kinetic energy metric $M$, i.e.
	\beq
	M^\sharp: T_p^*\mathcal{Q} \longrightarrow T_p \mathcal{Q} \quad \text{and} \quad M^\flat:T_p\mathcal{Q} \longrightarrow T_p^*\mathcal{Q}.
	\eeq
Similarly, we have the equivalent Riemannian problem
	\beq
	\label{geomcontrol1}
	\nabla_U U = G^\sharp\left[ u \right], \quad \text{with} \quad u=\sum_{i=1}^m u_i \hat \theta^{(i)} \quad \text{with} \quad \hat \theta^{(i)} = G^\flat \left[\hat e_{(i)} \right].
	\eeq
Note that $u$ is the same as in equation \eqref{controlEL} but expressed in terms of a different co-frame, the one corresponding to the Jacobi metric $G$. 

If  $\text{span}(\mathcal{D}) = T_p\mathcal{Q}$ then we say the system is \emph{fully-actuated}. Otherwise, we say it is \emph{under-actuated}. Most of the relevant situations in control theory involve under-actuated systems, i.e. those for which $\text{span}(\mathcal{D}) \subset T_p\mathcal{Q}$ \cite{auckly1999control}. Furthermore, if one provides an input force depending solely on time, controlling the evolution of the system in the configuration space $\mathcal{Q}$, then $u$ is called and \emph{open loop control}\footnote{In fact, open loop \emph{tracking} control for desired reference motion trajectories could be directly synthesized for a class of under-actuated controllable dynamical systems in absence of uncertainty, e.g. differentially flat systems \cite{fliess1995flatness}. In such case, system variables (states and control) can be expressed in terms of a set of flat output variables and a finite number of their time derivatives}.  However, the central interest in control theory is in those systems which can regulate themselves against unknown perturbations, that is, those whose control input is responsive to spontaneous variations of the system configuration. This is referred as \emph{closed loop control}.  In such case, the control input force depends on the state of the system at any given time, that is, $u$ must be a section of the co-tangent bundle
and the control objective is that the integral curves of \eqref{controlEL} -- or, equivalently, those  satisfying \eqref{geomcontrol1} -- remain ``close'' to a reference path $\gamma^* \subset \mathcal{Q}$ with some desired equilibrium properties.

To clearly state our Riemannian control problem, we propose a slight modification  of Lewis' work \cite{bullo2004geometric, Lewis:2007:WLD:1317140.1317151}. Moreover, without any loss of generality and to keep our argument sufficiently simple, we will restrict ourselves to the case where no external nor gyroscopic forces are present (cf. \cite{gharesifard2008geometric}). Thus, let us define an open loop control system as the triad 
	\beq
	\Sigma_{\text{ol}} \equiv \left\{\mathcal{Q},G_{\text{ol}},\mathcal{W} \right\},
	\eeq
where $\mathcal{Q}$ and $G_{\text{ol}}$ are the configuration space and the Jacobi metric, respectively;  and $\mathcal{W}\subset T^*\mathcal{Q}$ is  the control sub-bundle defined at each point as 
	\beq
	\mathcal{W}_p = \text{span}(F_p), \quad F_p = \left\{\hat \theta^{(i)}(p) \right\}_{i=1}^m \quad \forall p \in \mathcal{Q}, 
	\eeq
so that equation \eqref{geomcontrol1} above is satisfied for a certain $u \in \mathcal{W}$ given some reference $\gamma^*$. The goal is to obtain a pair -- the closed loop system -- 
	\beq
	\Sigma_{\text{cl}} = \left\{\mathcal{Q}, G_{\text{cl}} \right\},
	\eeq
such that the geodesics of $G_{\text{cl}}$  \emph{match}  the open loop solutions of \eqref{geomcontrol1}. Here, we use the term match to indicate that the geodesics of the closed loop metric are not required to coincide at every point with the solutions of the open loop system, but merely that the vector fields share the same singular points together with their equilibrium properties, i.e. that the equilibria of the open loop system are the same as those of the closed loop geodesic vector field and with the same asymtotic behaviour. In such case we say  that $u\in\mathcal{W}$ is  a stabilizing input for the geometry shaping problem solved by $G_{\text{cl}}$ \cite{Lewis:2007:WLD:1317140.1317151}.   Using this formulation the  stability properties of the closed loop system can be assessed directly by means of the Riemann tensor of $G_{\text{cl}}$ through the geodesic deviation equation. Intuitively, unstable regions should correspond to subsets of the configuration space where the eigenvalues of the Riemann tensor take negative values and the geodesics in a congruence are divergent, providing us with a completely geometric and coordinate independent  stability criterion.

Notice that the closed loop metric is not required to be a Jacobi metric, that is, it is not necessary to find a function $V_{\text{cl}} = V_{\text{cl}}(p)$ such that  $G_{\text{cl}}$ be of the form \eqref{JacobiMetric}. Moreover, it may not even share the same signature with $G_{\text{ol}}$ (cf. Remark 4.7 in \cite{gharesifard2008geometric}). Finally, recalling that the geodesics of the closed loop metric correspond to the integral curves of a purely kinetic Lagrangian vector field with kinetic energy metric $G_{\text{cl}}$, and that the open loop metric is a Jacobi metric, $G_{\text cl}$ must satisfy the partial differential equation
	\beq
	\label{pde1}
	G_{\text{ol}}^\sharp \left[\sum_{i=1}^m u_i(U_p) \hat \theta^{(i)}(p) \right] = \left(\nabla^{\text{cl}} - \nabla^{\text{ol}} \right)\left[U_p,U_p \right],
	\eeq 
where $\nabla^{\text{cl}}$ and $\nabla^{\text{ol}}$ are the Levi-Civita connections of $G_{\text{cl}}$ and $G_{\text{ol}}$, respectively;  we have used the shorthand  $\left(\nabla^{\text{cl}} - \nabla^{\text{ol}} \right)\left[U_p,U_p \right]$ to denote the \emph{connection difference tensor}, for every pair $(p,U_p)$ denoting a posible  state of the system at point $\tilde p \in T\mathcal{Q}$ [cf. equation \eqref{velocity}, above]. We have been emphatic on the fact that the stabilizing control input must depend on the state of the system. Equations of the form \eqref{pde1} might be used to define further geometric structures in the space they are defined, as has been done in the case of the geometrization programme of thermodynamics and fluctuation theory developed in \cite{bravetti2015sasakian, bravetti2015conformal, fernandez2015generalised}.
 
In recent years, solutions to equation \eqref{pde1} have been the object of various studies \cite{crasta2015matching,ng2013energy}. However, one should be aware that  the existence of a solution of \eref{pde1} satisfying some desired properties might be severely constrained by the topology of $\mathcal{Q}$ \cite{Lewis:2007:WLD:1317140.1317151}. Thus, in general,  there is no generic criterion for deciding the solvability of \eqref{pde1}. Nevertheless, it has been shown that the case of one degree of under-actuation is fully stabilizable by means of \eqref{pde1} \cite{gharesifard2011stabilization}; modulo the difference introduced by using the open loop Jacobi metric in the problem's definition. Interestingly, in such case, equation \eqref{pde1} resembles very closely that of the definition of the general relativistic \emph{elasticity difference tensor} introduced by Karlovini and Samuelsson \cite{karlovini2003elastic}
	\beq
	\label{edt}
	S = P(\nabla - \tilde \nabla)
	\eeq
where $\nabla$ and $\tilde \nabla$ represent the spacetime and the (\emph{pulled back}) matter space metrics, respectively (cf. \cite{carter1972foundations,lopez2011covariant,lopez2011thermal,andersson2007relativistic} for a detailed presentation of matter spaces in relativistic elasticity and dissipation) and $P$ denotes the orthogonal projection with respect to a  free falling geodesic congruence. Such a tensor has been used to recast the relativistic Euler's equations  in terms of the  Hadamard Elasticity, a form which is elegant and useful in the study of wave propagation in relativistic elastic media \cite{vaz2008analysing}. This provides us with an interesting link between relativistic elasticity and geometric control theory which deserves further exploration.

%Making the analogy between the control theoretic equation \eqref{pde1} and the elasticity difference tensor \eqref{edt} provides us with, on the one hand, a bridge between geometric control theory and relativistic elasticity, where new insights could arise from the techniques developed in \cite{GHARESIFARDsiam}; and, on the other, a visual interpretation of the open and closed loop Levi-Civita connections as characterising the ``elasticity'' associated with the integral curves of \eqref{geomcontrol1}.

Let us close this contribution by noting that there might be several stabilizing inputs for a given desiderata. In such case, one might look for the `most suitable' geometry solving our control objective and use \eqref{pde1} as a constraint for a variational problem. Its particular form  should be fixed by some \emph{a priori} known cost functional that is to be extremized by the sought for closed loop metric. Such is the standard \emph{optimal control problem} \cite{sethi2000optimal}. To preserve the geometric nature of the whole construction, the cost functional can only depend on scalars formed  from the metric and its derivatives, that is
	\beq
	\label{mfunc}
	\mathcal{A}\left[G \right] = \int_\mathcal{Q} \mathcal{F}\left(G,G',G'',... \right) \sqrt{G}\ \d^n q,
	\eeq 
where $\sqrt{G}\ \d^n q$ is the invariant volume element on $\mathcal{Q}$. Thus, our search for geometries solving a control objective has led us to the study of cost functionals akin to the various classes of gravitational theories where the extremum is achieved by geometries  required to be compatible with the \emph{observed} free fall motion of certain spacetime observers\cite{lovelock1971einstein}. Therefore,  a complete solution to our problem -- if it exists --, should be a metric $G_\text{cl}$ extremizing \eqref{mfunc}
%	\beq
%	\left.\frac{\delta \mathcal{A}}{\delta G} \right\vert_{G=G_{\text{cl}}} = 0
%	\eeq
such that equation \eqref{pde1} is satisfied for a given open loop control. An exploration of the equations of motion stemming from the class of cost functionals constructed from curvature invariants in control theory will be the subject of further investigations.
 In this sense, the lessons learned and the results obtained from the variational formulation of gravitational theories might find a novel application in the realm of geometric control theory.

\section*{Acknowledgements}

CSLM wishes to express his gratitude to the Organizing Committee for their kind hospitality and strong efforts in providing us with such an inspiring venue for this celebration. 

\section*{References}

%\bibliography{myrefs_control}{}
%\bibliographystyle{unsrt}

\end{document}